\documentclass[12pt]{amsart}
\usepackage{amsfonts} \usepackage{amsmath} 
\theoremstyle{plain}              
\newtheorem{Thm}{Theorem} \newtheorem{Lem}[Thm]{Lemma}      
\newtheorem{Prop}[Thm]{Proposition}                         \theoremstyle{definition}
                                
\theoremstyle{remark}     \newtheorem{Rem}{Remark}

\def\Z{\mathbb Z}

\begin{document}
\author{Selahi Durusoy} 
\title{Heegaard-Floer Homology and a family of Brieskorn spheres}
\address{Michigan State University}
\address{durusoyd@math.msu.edu}

\begin{abstract}We compute the Heegaard-Floer homology for the family $\Sigma(2,3,6n+1)$
of Brieskorn spheres using the algorithm given in \cite{os_plumbed}.
\end{abstract}

\maketitle \setlength{\unitlength}{.62mm} 

\setcounter{section}{-1}
\section{Introduction}

\noindent Ozsv\'ath and Szab\'o have given a combinatorial description for Heegaard-Floer homology in \cite{os_plumbed}.
Using this algorithm we compute Heegaard-Floer homology of $-\Sigma(2,3,6n+1)$:

\vspace{1mm}
{\Thm $HF^+(-\Sigma(2,3,6n+1)) = T^+_0 \oplus \Z^n_{(0)}$.}
\vspace{1mm}

\noindent This family of homology spheres can be obtained by doing $-1/n$ surgery on right handed trefoil knot. This family and some more were considered in \cite{fs} where their instanton Floer homology is calculated.
For the Heegaard-Floer homology computations for the family $\Sigma(2,2n+1,4n+3)$ see \cite{r}. Another
combinatorial description for calculating $HF^+$ has been given by N\'emethi in \cite{nem}.

\section{Remarks on the algorithm}

Consider a negative definite weighted graph $G$ with at most one bad vertex (i.e., $d(v)>-v \cdot v$ for at most one vertex $v$ where $d(v)$ denotes the number of edges containing the vertex $v$.)
Let $X(G)$ denote the plumbed disk bundle over  $S^2$ and $Y(G)$ its boundary.
In \cite{os_plumbed} the subset of $\mathbb{H}^+(-Y(G))$ of $\text{Map}(Char(G),T_0^+)$ of finitely supported maps
satisfying $U^{m+n}\cdot f(K+2PD[v]) = U^m \cdot f(K)$ whenever $min\{m,m+n \}\ge 0$ and $K\cdot v + v\cdot v = 2n$ 
is considered and the following is shown:
{\Thm \cite{os_plumbed} For such a graph $G$, for each $Spin^c$ structure $t$ over $-Y(G)$,
$$HF^+(-Y(G),t) \cong \mathbb{H}^+(G,t).$$}
\noindent We will be working with homology spheres, so we suppress the $Spin^c$ structure from the notation.
In the computations instead of working with elements of $\mathbb{H}^+(G)$, \cite{os_plumbed}
works with elements of $K^+(G)$, which is the set of equivalence classes of elements in $\Z^{\ge 0} \times Char(G)$,
with the equivalence relation defined by
$$(m,K) \sim (m+n,K+2PD[v])$$
where $v$ is a vertex in $G$ with $K\cdot v + v \cdot v = 2n$ and $min(m,m+n) \ge 0$. Equivalence class of
$(m,K)$ will be denoted by $U^m\otimes K$.

For an equivalence class $U^m\otimes K$, define its $U$-depth as the smallest number $l$ so that $(l,K')$
is a representative of $U^m\otimes K$ for some vector $K'$. $K^+(G)$ is determined by elements of $U$-depth $0$
and the $U$ action on $K^+(G)$, which follows from 

{\Prop (Prop 3.2 in \cite{os_plumbed}) For an equivalence class $U^m\otimes K$ of $U$-depth $0$, there is a unique
representative $(0,K)$ satisfying
\begin{equation}
v_i \cdot v_i +2 \le K \cdot v_i \le - v_i \cdot v_i \text{ for each } i.
\label{cond13} 
\end{equation}
\noindent Conversely if a vector $K$ satisfies (\ref{cond13}), then $K$ has $U$-depth $0$ if and only if
$K$ supports a good full path (in this case $K$ will be called a basic vector).}

\noindent In the above, full path stands for a path of vectors $K_1, K_2, ..., K_n$ in $Char(G)$ with $K_1$ satisfying (\ref{cond13}) obtained by adding $2PD[v_i]$ if $K_i\cdot v_j + v_j\cdot v_j = 0$ for some $j$, until 
$-K_n$ satisfies (\ref{cond13}) or $K_n\cdot v_j + v_j\cdot v_j > 0$ for some $j$. It is called good if  $-K_n$ satisfies (\ref{cond13}).

In the proof of \cite{os_plumbed}, Proposition 3.2, it is shown that given a characteristic vector
$M$, the final vector of any full path is identical, hence we observe the following useful

{\Rem \label{hereditary} Observe that if a vector supports a good full path, then all full paths are good,
hence finding one bad full path means the initial vector is not basic. Secondly observe that
bad full paths are hereditary, hence if a vector for a subgraph has a bad full path, then so does the
containing vector and graph.}

This will reduce the number of vectors that we need to check if they support a good full path or not.
We will express the vectors $K \in Char(G)$ as sequences $(K \cdot v_i)$.

{\Lem \label{a_s} For the linear graph $A_s$ with $s$ vertices and each weight $-2$, there are no good full paths starting at vectors $K$ satisfying (\ref{cond13}) with $K \cdot v_i = 2$ for more than one $i$.}

\noindent {\em Proof.} We will use induction on $s$, observing that we can use hereditary property of bad full paths.
(\ref{cond13}) implies $K \cdot v_i \in \{0,2\}$. For $s=2$, $(2,2)$ for $A_2$ has a bad full path obtained by adding $2PD(v_1)$. For $s>2$, observe that if 
$K \cdot v_i = K \cdot v_j =2$ for some $i \ne j$, then $K$ is equivalent
to a vector containing $(2, 2)$ as shown below, hence has a bad full path.
$$( *, 2, 0, 0, ..., 0, 0, 2, *') \sim ( *'', -2, 2, 0, ..., 0, 0, 2, *') \sim ( *'', -2, 0,..., -2, 2,2,*'). \square$$

\section{The family $\Sigma(2,3,6n+1)$}

Consider the family of Brieskorn spheres $Y(n) = \Sigma(2,3,6n+1)$.
The negative definite plumbing graph defining $Y(n)$ is a tree with weights $-1$ on central node, $-2, -3, -7$
on adjacent nodes and a $-2$ chain of length $n-1$ starting at $-7$ as follows:
%
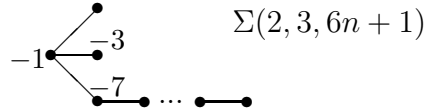
\begin{center}
\begin{picture}(100,30)(0,0) \put(7,12){\circle*{2}} \put(2,11){\makebox(0,0){{\small $-1$}}}
 \put(7,12){\line(1,1){10}}  \put(17,22){\circle*{2}}
\put(7,12){\line(1,0){10}}  \put(17,12){\circle*{2}}
\put(19,15){\makebox(0,0){{\small $-3$}}}
 \put(7,12){\line(1,-1){10}}  \put(17,2){\circle*{2}}
 \put(17,2){\line(1,0){10} \circle*{2}}
 \put(39,2){\circle*{2}}
 \put(39,2){\line(1,0){10} \circle*{2}}
 \put(19,5){\makebox(0,0){{\small $-7$}}}
\put(46,18.5){\makebox(0,0)[l]{$\Sigma(2,3,6n+1)$}}
\put(30,2){\makebox(0,0)[l]{$...$}}
\end{picture}
\end{center}

The Heegaard-Floer homology of the first member of this family, $HF^+(-\Sigma(2,3,7))$, has been studied in \cite{os_plumbed}, \cite{r}.

{\Lem For arbitrary $n$, the basic vectors for $Y(n) = \Sigma(2,3,6n+1)$ are
\begin{eqnarray}
K_1 & = & (1, 0, -1, -5, 0, 0, 0, ... , 0)\nonumber \\
K_2 & = & (1, 0, -1, -3, 0, 0, 0, ... , 0)\nonumber \\
K_3 & = & (1, 0, -1, -5, 2, 0, 0, ... , 0)\nonumber \\
K_4 & = & (1, 0, -1, -5, 0, 2, 0, ... , 0)\nonumber \\
K_5 & = & (1, 0, -1, -5, 0, 0, 2, ... , 0)\nonumber \\
    & \vdots & \nonumber \\
K_{n+1} & = & (1, 0, -1, -5, 0, 0, 0, ... , 2) \nonumber 
\end{eqnarray}
} 

\noindent {\it Proof.} Clearly for each $j$, $K_j$ satisfies (\ref{cond13}).
Next we need to see that among all characteristic vectors satisfying (\ref{cond13})
these are the only ones supporting good full paths.
For $n=1$ this is done in \cite{os_plumbed}, and we verified it by computer.

By remark \ref{hereditary}, for $n>1$ first $4$ entries of
a basic vector has to coincide with one of $(1, 0 -1, -3), \, (1, 0, -1, -5)$ which were
computed in \cite{os_plumbed} for $n=1$. Other entries are either $0$ or $2$.
Moreover lemma \ref{a_s} implies that for basic vectors $K$ for $Y(n)$ there can be
at most one vertex with $K \cdot v_i = 2$.

\noindent {\it Claim.} $(1, 0, -1, -3, * )$ has a bad path if $*$ has a non-zero entry.

\noindent {\it Proof of claim.} As in lemma \ref{a_s}, we can find a vector $(1, 0, -1, -3, 2, *')$ equivalent to $K$.
But $(1, 0, -1, -3, 2)$ has a bad path obtained by addind $2PD(v_i)$ in the order
\mbox{$i= 1, 2, 1, 3, 1, 2, 1, 5, 4, 1, 2, 1$} and  bad paths are hereditary. $\square$ 

\noindent Therefore the only basic vector with initial segment $(1, 0, -1, -3)$ is $(1, 0, -1, -3, 0, ..., 0)$.

Next we need to show that $K_1, ..., K_{n+1}$ support good paths. This we do by explicitly giving the paths.
First, $1, 2, 1, 3, 1, 2, 1$ is a good full path for both $K_1$ and $K_2$. For others, the path starts the same, but continues as:

\begin{eqnarray}
5, \;\; 6,  \;\; 7, ..., n+3 & \text{ for } K_3 \nonumber \\
6,5, \;\; 7,6, \;\; 8,7, ..., n+3, n+2 & \text{ for } K_4 \nonumber \\
7,6,5, \;\; 8,7,6, \;\; 9,8,7, ..., n+3, n+2, n+1 & \text{ for } K_5 \nonumber \\
\vdots & \nonumber \\
n+3, n+2, n+1, ..., 5 & \text{ for } K_{n+1} \nonumber 
\end{eqnarray}

\noindent This finishes the proof of the lemma. $\square$

For each $K_i$, when we compute the renormalized lengths $\frac{K \cdot K + |G|}{4}$, each time we get $0$.
Next we investigate relationships between $U$ powers of $K_i$.

{\Lem $U\otimes K_i \sim U \otimes K_j \sim L = (-3, 2, 5, 1, 0, 0, ..., 0)$ for $1 \le i,j \le n+1$ }

\noindent {\it Proof.} For $K_1$, the sequence $1,1,2,1$ leads to $L$.

\noindent For $i>1$, the path from $K_i$ leading to $L$ is of the following form:
$$ 1,1,2,1,2,3,1, \;\; A_{n,i}, \;\; B_n$$
where $B_n$ is $1,2,3,1,4,1, 2, 1$ followed by $5,6,...,n+3$ followed by $1$ and $A_{n,i}$ is of the following form:

$$C_{n-i+1}, C_{n-i+2}, ..., C_{0}.$$

\noindent In the above, $C_k$ denotes the sequence $1,2,3,1,4,1,2,1,3,1,2,1, \;\; 5,6,4+k$ if $k>0$ and empty path if $k=0$. As an example, for $n=4$, the path from $K_2$ to $L$ is given by

\begin{eqnarray}
1,1,2,1,2,3,1, A_{4,2}, B_4 & = &1,1,2,1,2,3,1, C_3, C_2, C_1, B_4 \nonumber \\ 
& = & 1,1,2,1,2,3,1, \nonumber \\
& & 1,2,3,1,4,1,1,2,1,2,3,1,  \nonumber \\
& & 1,2,3,1,4,1,2,1,3,1,2,1,5, \nonumber \\
& & 1,2,3,1,4,1,2,1,3,1,2,1,5,6, \nonumber \\
& & 1,2,3,1,4,1,2,1,5,6,7,1 \nonumber 
\end{eqnarray}

\noindent It is straightforward to check that these paths end at $L$. Next one observes that
as in \cite{ad}, throughout these paths the $U$-depth stays between $0$ and $1$, hence we get
$U\otimes K_i \sim L$ as announced. $\square$

Now we know that $K^+(G)$ consists of $U^0\otimes K_1, U^0\otimes K_2$ and $U^m\otimes K_1$ for $m>0$.
For any $f$ in $\mathbb{H}^+(G)$, $f(K_1)\in T_0^+$ determines the images $\tilde{f}(m,K_1)$ of the induced map $\tilde{f}:\Z \times Char(G) \rightarrow T^+_0$.
The remaining values $f(K_i)$ for $i>1$ are also determined up to addition of an element of $\Z_{(0)}$. This finishes the
proof of the theorem.
 $\square$

\vspace{2mm}

\noindent {\em Acknowledgements.} I would like to thank Selman Akbulut for advice, support and encouragement.
Also I wish to thank John McCarthy and Yildiray Ozan for discussions and Dina Shnaider for encouragement.

 \end{document}